\documentclass{article}%
\usepackage{amsfonts}
\usepackage{amsmath}
\usepackage{amsthm}
\usepackage{amssymb}
\usepackage{supertabular}%
\usepackage{xcolor}
\def\vep{{\varepsilon}}

\def\cP{\mathcal{P}}

\def\Q{{\mathbb Q}}
\def\Z{{\mathbb Z}}

\def\bI{{\bf I}}

\def\al{{\alpha}}
\def\be{{\beta}}
\def\ga{{\gamma}}
\def\de{{\delta}}
\def\la{{\lambda}}
\def\La{{\Lambda}}
\def\z{{\zeta}}
\def\th{{\theta}}
\def\ep{{\epsilon}}
\def\vep{{\varepsilon}}

\def\si{{\sigma}}
\def\om{{\omega}}
\def\Ht{{\mathrm{h}}}

\def\val{\mathrm{v}}
\def\wal{\mathrm{w}}

\def\h{{\rm h}\kern.5pt}

\def\gcd{{\mathrm{gcd}}}
\def\res{{\mathrm{Res}}}
\def\dis{{\mathrm{Disc}}}

\def\proof{{\bf Proof}.\:}
\def\proofend{\hspace*{1mm} \hfill{$\Box$}}

\newcommand{\bincf}[2]{\left(\!\!\!\begin{array}{c} #1 \\ #2 \end{array}\!\!\!\right)}
\newcommand{\ideal}[1]{\langle #1 \rangle}

\newcommand\jac[2]{\left(\frac{#1}{#2}\right)}

\newcommand{\omfroot}[1]{\frac{\om}{p}\cdot\frac{1+\z^{#1}}{1-\z^{#1}}}

\newtheorem{theorem}{Theorem}[section]

\newtheorem*{conjecture}{Conjecture}
\newtheorem{corollary}[theorem]{Corollary}

\newtheorem{lemma}[theorem]{Lemma}

\newtheorem{proposition}[theorem]{Proposition}

\begin{document}
\title{On $px^2 + q^{2n}= y^p$ and related Diophantine equations}
\author{A.~Laradji\thanks{Department of Mathematics \& Statistics, KFUPM, Dhahran,
Saudi Arabia, e-mail: alaradji@kfupm.edu.sa} 
\and
M.~Mignotte\thanks{Universit\'{e} Louis Pasteur,U.~F.~R. de Math\'{e}matiques, Strasbourg, France,
e-mail: mignotte@math.u-strasbg.fr}
\and
N.~Tzanakis\thanks{Department of Mathematics, University of Crete, Iraklion,
Greece, e-mail: tzanakis@math.uoc.gr, http://www.math.uoc.gr/\~{}tzanakis}
}
\date{{\small Revised version, October 8, 2020}\,\footnote{Corrected and/or added text
in \textcolor{blue}{blue}.}}

\maketitle
\begin{abstract}
The title equation, where $p>3$ is a prime number 
$\not\equiv 7 \pmod 8$, $q$ is an odd prime number and $x,y,n$ are positive integers with $x,y$ relatively prime, 
is studied. When $p\equiv 3\pmod 8$, we prove (Theorem 2.3) that there are no solutions. 
For $p\not\equiv 3\pmod 8$ the treatment of the equation turns out to be a difficult task. 
We focus our attention to $p=5$, by reason of an article by F.~Abu Muriefah, published in this journal, vol.~128 
(2008), 1670-1675. Our main result concerning this special equation is Theorem 1.1, whose proof is based on 
results around the Diophantine equation $5x^2-4=y^n$ (integer solutions), interesting in themselves, 
which are exposed in Sections 3 and 4. These last results are obtained by using tools such as Linear Forms 
in Two Logarithms and Hypergeometric Series.
\end{abstract}
\section{Introduction} \label{intro}
Diophantine equations of the form $px^2+c=y^p$,  where $c$ is a nonzero integer and $p$ is an odd prime, 
have been studied by several authors. When $c=2^n$, the case  $p=3$ was solved by Rabinowitz in \cite{RAB}, 
while Le dealt with the case $p>3$ in \cite{LE}. The case $c=3^n$ was considered by Abu Muriefah in \cite{AMT}. 
Cao \cite{CAO} treated the cases $c=1$ and $c=4^n$ (see also \cite{AMI}, \cite{ADA}, \cite{BROWN}, \cite{COHN} 
for closely related results). We should moreover mention that the equation has no solution in positive integers 
$x,y$ when $c=-1$, as can be inferred from the work of Nagell \cite{NAG} and Cao \cite{CAO1}.

The case when $c=q^{2n}$,  where $q$ is an odd prime, was studied in the recent paper \cite{AM}
of Abu Muriefah. 
Let us first note that, for fixed $n$ and $p\geq 5$, the Diophantine equation $px^2+q^{2n}=y^p$ 
and, more generally, the Diophantine equation $pX^2+Y^{2n}=Z^p$, have at most finitely many solutions $(x,q,y)$ 
and $(X,Y,Z)$, respectively. Indeed, in this case, $\frac{1}{2}+\frac{1}{2n}+\frac{1}{p}<1$ and the claim
follows from Theorem 2 of \cite{DaGra}. 
A main result in the aforementioned paper \cite{AM}, namely, Theorem 3.1, states that the equation 
$5x^2+q^{2n}=y^5$, has two families of solutions given by 
$y=\phi_{3k},\phi_{3k+1}$ (or $\psi_{3k+1},\psi_{3k+2})$, \ $k>1$,  
where $\phi_{k}$ (respectively $\psi_{k}$) is the $k$th term of the Fibonacci
sequence (respectively the Lucas sequence). However, straightforward
computations show that the only Fibonacci or Lucas number $y<1000$ satisfying
the title equation when $p=5$ is $y=21$ with $x=410$,  $q=1801$,  $n=1$,  and, further, in
$5\times183630^{2}+160201^{2}=181^{5}$, $160201$ is prime and $181$ is neither
a Fibonacci nor a Lucas number. The same Theorem 3.1 of \cite{AM} states that, if $p>3$ is a prime
$\not\equiv 7\pmod 8$ and $q$ is another odd prime, then there are no integer solutions $(x,y,n)$ 
to the equation $px^2+q^{2n}=y^p$ with $(x,y)=1$. The proof of Theorem 3.1 in \cite{AM}, just before its end,
contains an obvious, non-rectifiable error, at case 2, when $-16apb^2$ is set equal to $-16apq^{2m}$
although $b=\pm q^j$ with $0\leq j<m$. That the said proof is erroneous is also pointed out by P.G.~Walsh in
his review \cite{Walsh} of \cite{AM}. One of our aims in the present paper is to prove Abu Muriefah's assertion 
when $p\equiv 3\pmod 8$; see Theorem \ref{th p equiv 3 mod 8}. Our proof, rather than rectifying Abu Muriefah's 
argument (this is probably impossible), goes through totally different lines. Unfortunately, our arguments
cannot be extended to the case $p\equiv 1\pmod 4$. 

As we revisited the title equation, we further discovered some new results, like Theorems \ref{th on main eq}, 
\ref{th gen exp eq} and \ref{th with y=q} which, we believe, merit one's attention. Moreover, since the
powerful techniques of sections \ref{section exp eq with y prime} and \ref{section general exp eq} 
are also applicable (after the appropriate modifications) to Diophantine equations other than the ones treated 
in this paper, we thought it useful to expose them in some detail, enough for the reader to profit from them.

As we stated above, one of the main results of this paper is the following
\begin{theorem} \label{th on main eq}
 Let $q$ be an odd prime. If either $q\not\equiv 1\pmod{600}$ or $q \leq 3\cdot 10^9$, then there is
no integer solution $(x,y,n)$ to the equation
\begin{equation} \label{eq initial}
 5x^2+q^{2n} = y^5\,,\quad x,y,n >0\,.
\end{equation}
Otherwise, there exists at most one integer solution $(x,y,n)$ and if it actually exists, then 
it must satisfy the following conditions:
\\
(i) $n<820$ and $\gcd(n,2\cdot 3\cdot 5\cdot 7\cdot 11\cdot 13)=1$.
\\
(ii) There exists an integer $v$ such that $x=10v(80v^{4}-40v^{2}+1)$,
$ y=20v^{2}+1$, $ q^{n}=2000v^{4}-200v^{2}+1$.
\end{theorem}
{\em Remark}. If the prime $q$ is of the form $q=2000v^4-200v^2+1$ (the first few primes of this shape
are 1801, 160201, 1245001, 4792201, 8179201), then $(x,y,n)=(10v(80v^4-40v^2+1),20v^2+1,1)$ is a solution
to (\ref{eq initial}) and, according to the theorem, this is the only one (with $x,y>0$). We have not been
able, however, to find a prime $q$ such that the corresponding equation (\ref{eq initial}) has a solution
with $n>1$. Therefore, we state the following
\begin{conjecture}
If the prime $q$ is not of the form $q=2000v^4-200v^2+1$, then the equation (\ref{eq initial}) has
no solutions.
\end{conjecture}
The proof of Theorem \ref{th on main eq} follows from a straightforward combination of 
Corollary \ref{proposition on main eq} and Theorem \ref{th with y=q}, a second main result of our paper.
In turn, the proof of Theorem \ref{th with y=q} relies on a third main result, namely, Theorem \ref{th gen exp eq}
concerning the equation $5x^2-4=y^n$ which is interesting for its own sake. 
Indeed, in recent years, important papers are devoted to equations of the form $x^2+C=y^n$.
One main strategy for attacking such equations is based on the so called {\em modular method} which has been 
successfully applied in quite a number of cases; see Chapter 15 (by S.~Siksek) in H.~Cohen's book \cite{Cohen},
the survey article \cite{FSAbu} and \cite{BMSII} and the references therein. 
For our equation $5x^2-4=y^n$, the existence of the trivial solution $(x,y)=(1,1)$ makes the modular approach
unsuccessful and prevents us from giving the complete solution $(x,y,n)$. Thus, our Theorem \ref{th gen exp eq} 
offers only a partial result which, at present, seems to be best possible. 
\section{The Diophantine equation $px^2+q^{2n}=y^p$}
The main results of this section are Proposition \ref{proposition general Muriefah} and 
Theorem \ref{th p equiv 3 mod 8}
\begin{proposition} \label{proposition general Muriefah}
Let $p>3$ be a prime number $\not\equiv 7 \pmod 8$ and let $q$ be an odd prime number.
If $x,y,n$ are positive integers with $x,y$ relatively prime such that
\begin{equation} \label{eq 3.1}
px^2+q^{2n}=y^p\,, \quad (x,y)=1, \ n>0\,,
\end{equation}
then there exists a rational integer $a$ such that
\begin{equation} \label{q^n}
 \pm q^n = \sum_{i=0}^{(p-1)/2} \bincf{p}{2i+1}a^{p-2i-1}(-p)^{(p-2i-1)/2} 
\end{equation}
and
\begin{equation} \label{x}
 x = \sum_{i=0}^{(p-1)/2} \bincf{p}{2i}a^{p-2i}(-p)^{(p-2i-1)/2}
\end{equation}
\end{proposition}
\proof
The condition $(x,y)=1$ implies $p\neq q$ and the condition $p\not\equiv 7\pmod 8$ implies 
that $y$ is odd.

We work in the imaginary quadratic field $K=\Q(\om)$, where $\om=\sqrt{-p}$. Equation (\ref{eq 3.1}) factorizes
as $(\om x+q^n)(-\om x+q^n)=y^p$ and, trivially, the factors in the left-hand side are relatively prime.
This implies an ideal equation $\ideal{\om x+q^n}={\bf I}^p$, where ${\bf I}$ is an integral ideal of $K$.
Since the ideal-class number of $K$ is strictly less than $p$ (see page 199 of \cite{Fai}), the 
above ideal equation implies that ${\bf I}$ is a principal ideal, therefore we obtain the equation
\begin{equation} \label{xomega+q^n}
x\om +q^n = \al ^p\,,\quad \al=a\th+b\,,
\quad \th =\begin{cases}
               \om & \mbox{if $p\equiv 1,5\pmod 8$} \\
               \frac{1+\om}{2} & \mbox{if $p\equiv 3\pmod 8$}
              \end{cases}
\end{equation}
for some rational integers $a$ and $b$.
In case $p\equiv 3\pmod 8$ we write the above equation as $q^n-x +2x\th=(b+a\th)^p$, implying
$(b+a\th)^p\equiv\mbox{$0$ or $1\pmod 2$}$. From this, we see that $a$ cannot be odd. For, otherwise, 
we would have $\th^p$ or $(1+\th)^p\equiv\mbox{$0$ or $1\pmod 2$}$. But we easily check that, for 
$k\not\equiv 0\pmod 3$, it is true that $\th^k, (1+\th)^k\equiv \mbox{$\th$ or $1+\th\pmod 2$}$, a contradiction.
Therefore, $a$ in (\ref{xomega+q^n}) is even and (\ref{xomega+q^n}) is equivalent to the simpler
equation
\begin{equation} \label{xomega+q^n simpler}
x\om +q^n = \al ^p\,,\quad \al=a\om+b\,,
\end{equation}
for some rational integers $a$ and $b$ of opposite parity since $y=a^2p+b^2$ is odd. Also, it is easy to see 
that $(pa,b)=1$. If we put 
\[
\be=-\bar \alpha\,,
\]
then
\[
x = \frac{\al^p + \be^p}{2\om}\,, \quad \frac{q^n}{b}=\frac{\al^p-\be^p}{\al-\be}
\]
and the fact that $(\al^p-\be^p)/(\al-\be)$ is an algebraic integer implies that $b$
divides $q^n$ (in $\Z$), hence $b=\pm q^j$ for some $j\in\{0,1,\ldots,n\}$.

At this point we note that the pair $(\al,-\be)$ is a Lehmer pair for which
\[
(\al^2-\be^2)^2=-16pa^2b^2\,.
\]
Concerning $j$ appearing in the relation $b=\pm q^j$ (see a few lines above), we distinguish two cases.
\\
(\rm{i}) $j>0$. Then, in the terminology of \cite{BHV}, $(\al,-\be)$ is a $p$-defective pair. By Theorems 1.4
and C of \cite{BHV} it easily follows that $p=5$ is the only possibility. Then, by Theorem 1.3 of \cite{BHV},
either $20a^2=\phi_{k-2\ep}$ for some $k\geq 3$, or $20a^2=\psi_{k-2\ep}$ for some $k\neq 1$, where
$\ep\in\{-1,1\}$, $(\phi_n)_{n\geq 0}$ denotes the Fibonacci sequence and $(\psi_n)_{n\geq 0}$
is defined by $\psi_0=2$, $\psi_1=1$ and $\psi_n=\psi_{n-1}+\psi_{n-2}$ for $n\geq 2$. It is easily checked that,
for every $n\geq 0$, $\psi_n\not\equiv 0\pmod 5$; therefore the second alternative must be excluded.
On the other hand, by Th\'{e}or\`{e}me 1.3 of \cite{BMSIII}, a relation of the form $\phi_k=5z^m$ with $m>1$
and $z>1$ is impossible, which excludes the first alternative as well.
\\
(\rm{ii}) $j=0$, so that $b=\pm 1$. Then, equating rational and irrational parts in (\ref{xomega+q^n simpler}), 
we respectively obtain the relations (\ref{q^n}) and (\ref{x}).
\proofend 
\begin{corollary} \label{proposition on main eq}
If the integers $x,y,q,n$,  where $(x,y)=1$,  $q$ is an odd prime and $n\geq1$,  satisfy the equation
\begin{equation} \label{Abu eq}
5x^{2}+q^{2n}=y^{5}
\end{equation}
then $(n,6)=1$,  $q\equiv 1\pmod{600}$ and there exists an integer $v$ such that
\begin{equation}\label{shape of sols}
x=10v(80v^{4}-40v^{2}+1)\,,\quad y=20v^{2}+1\,,\quad q^{n}=2000v^{4}-200v^{2}+1\,.
\end{equation}
\end{corollary}
\proof Applying Proposition \ref{proposition general Muriefah} with $p=5$ we obtain
$\pm q^n=125a^4-50a^2+1$ and $x=5a(5a^4-10a^2+1)$. Obviously, the minus sign in the first equation
is rejected \textcolor{blue}{and $a$} is even. Putting $a=2v$ in these relations we obtain the first and third
relation in (\ref{shape of sols}), and then the second relation results immediately. 
\\
We claim that $n$ is odd. Indeed, otherwise $(a,q^{n/2})$ would be an integral point on the elliptic curve 
defined by $Y^2=125X^4-50X^2+1$. But this elliptic curve has zero rank and its only rational point is 
$(X,Y)=(0,\pm 1)$, which forces $q=1$, a contradiction. 
\\
We also claim that $n$ is prime to 3. Indeed, let us write the third equation (\ref{shape of sols})
as $q^n+4=5w^2$, where $w=20v^2-1$. If $n$ were divisible by 3, then the last equation could be written as
$(5^2w)^2=(5q^{n/3})^3+500$, again forcing $q=1$, because it is well known since long
(see, for example, Table 8 in \cite{LoFi}) that the only integral solutions $(X,Y)$ to $Y^2=X^3+500$ are
$(X,Y)=(5,\pm 25)$.
\\
Finally, we show that $q\equiv 1\pmod{600}$. First, we write the third equation (\ref{shape of sols})
as $q^n=200v^2(v^2-1)+1$, which shows that $q^n\equiv 1\pmod{600}$. 
Let $r$ be the order of $q$ modulo $600$. Then $r$ divides $n$,  and since $\varphi(600)=160$ and $n$ is odd, 
we obtain $r=1$ or $r=5$. The latter case cannot hold, for, otherwise, 
$5x^{2}=y^{5}-(q^{2n/5})^{5}$,  which has no proper solutions by Th\'{e}or\`{e}me 2(2) of \cite{Kraus}.
We therefore conclude that $r=1$ and $q\equiv 1\pmod{600}$, as claimed.
\proofend
%
\begin{theorem} \label{th p equiv 3 mod 8}
Let $p,q$ be odd primes with $p\equiv 3\pmod 8$ and $p>3$. Then, the Diophantine equation
\begin{equation} \label{eq p equiv 3 mod 8}
 px^2 +q^{2n} = y^p
\end{equation}
has no positive integer solutions $(x,y,n)$ with $(x,y)=1$.
\end{theorem}
\proof By the relation (\ref{q^n}) we see that the equation (\ref{eq p equiv 3 mod 8}) implies
\[
 \pm q^n = \frac{(a\om +1)^p -(a\om-1)^p}{2}\,.
\]
Let us consider the polynomial
\[
 f(x):=\frac{(\om x +1)^p -(\om x-1)^p}{2}\,.
\]
Clearly, this is a polynomial in $\Z[x]$ of degree $p-1$, with leading coefficient $-p^{(p+1)/2}$ 
and constant term $1$.
\begin{quote}
{\em First Claim}: The polynomial $f(x)$ factorizes over $\Q[x]$ into two
relatively prime polynomials $f_1(x),f_2(x)\in\Z[x]$, each of degree $(p-1)/2$.
\end{quote}
Proof: Let $\z$ be a primitive $p$-th root of unity, i.e. a root of the $p$-th cyclotomic polynomial
$\Phi(x)=x^{p-1}+\cdots+x+1 = (x^p-1)/(x-1)$. Let also $g$ be a primitive root $\bmod\,p$.
Observe that the field $L=\Q(\z)$ contains $\om$. Indeed, it is a well-known fact that the Gauss sum
$G=\sum_{t=1}^{p-1}\jac{t}{p}\z^t$ satisfies $G^2=(-1)^{(p-1)/2}p$. Therefore, we can assume $\om=G$.
Below we will use the well-known fact that, in the field $L$ we have the factorization into prime ideals 
$\ideal{p}=\ideal{\la}^{p-1}$, where $\la=1-\z$. 
For $\be\in L$ we will write $\wal(\be)$ to denote $\val_{\la}(\be)$, the exponent of $\la$ in the
prime factorization of $\be$; and for $b\in\Q$ we will write $\val(b)$ to denote $\val_{p}(b)$, 
the exponent of $p$ in the prime factorization of $b$.
\\
We have
\[
 f(x)=\prod_{k=1}^{p-1}\{(\om x+1)-\z^k(\om x-1)\}\,,
\]
from which it follows that the roots of $f(x)$ are exactly the following:
\[
 \frac{\om}{p}\cdot\frac{1+\z^k}{1-\z^k}\,,\quad k=1,\ldots,p-1\,.
\]
Therefore,
\[
 f(x)=-p^{(p+1)/2}\prod_{k=1}^{p-1}\left(x- \frac{\om}{p}\cdot\frac{1+\z^k}{1-\z^k}\right)\,.
\]
Let us put now
\begin{eqnarray*}
 f_1(x) & = & 
p^{(p+1)/4}\prod_{j=0}^{(p-3)/2}\left(x-\frac{\om}{p}\cdot\frac{1+\z^{g^{2j}}}{1-\z^{g^{2j}}}\right)
\\
f_2(x) & = & 
-p^{(p+1)/4}\prod_{j=0}^{(p-3)/2}\left(x-\frac{\om}{p}\cdot\frac{1+\z^{g^{2j+1}}}{1-\z^{g^{2j+1}}}\right)
\end{eqnarray*}
so that $f(x)=f_1(x)f_2(x)$. We now show that the polynomials $f_i(x)$ have rational coefficients.
\\
The Galois group of the extension $\Q(\z)/\Q$ is cyclic generated by the automorphism $\si$, defined by
$\si(\z)=\z^g$. Since $\om\in\Q(\z)\setminus\Q$, we must have $\si(\om)\neq\om$, therefore, $\si(\om)=-\om$.
Consequently, for the typical root of $f_1(x)$ we have
\begin{equation} \label{action of sigma}
 \si\left(\frac{\om}{p}\cdot\frac{1+\z^{g^{2j}}}{1-\z^{g^{2j}}}\right) =
-\frac{\om}{p}\cdot\frac{1+\z^{g^{2j+1}}}{1-\z^{g^{2j+1}}} =
\frac{\om}{p}\cdot\frac{1+\z^{g^{2j'}}}{1-\z^{g^{2j'}}}\,,
\end{equation}
where
\[
 j'\equiv \frac{p+1}{4}+j\pmod{\frac{p-1}{2}}
\]
and we choose
\[
 j'=\begin{cases}
      \frac{p+1}{4}+j & \mbox{if $0\leq j < \frac{p-3}{4}$}  \\
      \frac{p+1}{4}+j-\frac{p-1}{2} & \mbox{if $\frac{p-3}{4}\leq j\leq\frac{p-3}{2}$}
    \end{cases}\,,
\]
so that $j'$ runs (exactly once) through all values $0,1,\ldots,\frac{p-3}{2}$ as $j$ runs through
these values. Consequently, the coefficients of the polynomial $f_1(x)$ are fixed by $\si$, which implies
that they belong to $\Q$; and similarly for $f_2(x)$. Actually, the coefficients of $f_1(x),f_2(x)$ are 
integers and we prove this as follows.
\\
First, we show that the absolute value of the constant coefficient of both $f_1(x)$ and $f_2(x)$ is 1.
Indeed, let $b_i$ be the constant coefficient of $f_i(x)$. We already know that $b_i\in\Q$. 
Moreover, multiplying the right equalities (\ref{action of sigma}) for $j=0,\ldots,(p-3)/2$ and then the
resulting products in the two sides by $-p^{(p+1)/4}$, we obtain $b_2=b_1$.
But, $b_1b_2$ is equal to the constant term of $f(x)$, which is $1$. Therefore $1=b_1b_2=b_1^2$, 
from which $b_1=b_2=\pm 1$.
\\
Let us put now
\[
 g(x)=x^{p-1}f(\frac{1}{x})\,,\quad  g_1(x)=x^{(p-1)/2}f_1(\frac{1}{x})\,,\quad
                            g_2(x)=x^{(p-1)/2}f_2(\frac{1}{x})\,,
\]
i.e.~these are the reciprocal polynomials of $f(x)$, $f_1(x)$ and $f_2(x)$, respectively.
Since $f(x)=f_1(x)f_2(x)$, we also have $g(x)=g_1(x)g_2(x)$.
Since the constant term of $f(x)$ is $1$, $g(x)$ has leading coefficient $1$; and since $f(x)$
has integer coefficients, so does $g(x)$. Therefore, the roots of $g(x)$ are algebraic integers. 
Analogously, the polynomials $g_i(x)$ have leading coefficients equal to $\pm 1$, their 
coefficients are rational numbers and their roots, being roots of $g(x)$, are algebraic integers. 
Therefore, these polynomials have coefficients in $\Z$; consequently, the same is true for the polynomials
$f_i(x)$, as claimed. 
\\
Now, observe that $f_1(x)$ and $f_2(x)$ have no common roots, therefore they are relatively prime.
\begin{quote}
{\em Second Claim}: Let $\res(f_1,f_2)$ be the resultant of the polynomials $f_1,f_2$. Then
\begin{equation}  \label{resultant f1,f2}
 \res(f_1,f_2)=\pm 2^{(p-1)^2/4} p^{(p^2-1)/8}\,.
\end{equation}
\end{quote}
Proof: We use the symbol $\dis$ to denote the discriminant. We have
\[
 \dis(f)=\dis(f_1f_2)=\dis(f_1)\dis(f_2)\res(f_1,f_2)^2\,.
\]
By the right-most equality in (\ref{action of sigma}) and the comments following it we see that 
$f_2(x)=0$ iff $f_1(-x)=0$, hence
\begin{equation} \label{disc and res}
 \dis(f)=\dis(f_1f_2)=\dis(f_1)^2\res(f_1,f_2)^2\,.
\end{equation}
Calculation of $\dis(f)$:
\begin{align*}
 \dis(f) & = p^{(2p-4)(p+1)/2}\prod_{1\leq i<j\leq p-1}\left(\omfroot{i}-\omfroot{j}\right)^2 \\
         & = p^{(p+1)(p-2)}\left(\frac{-1}{p}\right)^{(p-2)(p-1)/2}
             \prod_{1\leq i<j\leq p-1}\left(\frac{2(\z^i-\z^j)}{(1-\z^i)(1-\z^j)}\right)^2 \\
& = -2^{(p-1)(p-2)}p^{(p-2)(p+3)/2}\prod_{1\leq i<j\leq p-1}\left(\frac{\z^i-\z^j}{(1-\z^i)(1-\z^j)}\right)^2\,.
\end{align*}
The right-most product in the last equality is a unit times $\la^{-(p-1)(p-2)}$, therefore, 
$\wal(\dis(f))=(p-1)(p-2)(p+3)/2 - (p-1)(p-2)= (p-2)(p^2-1)/2$.
But, since $\dis(f)$ is a rational integer, it follows that
\begin{equation} \label{discr f}
 \dis(f)=\pm 2^{(p-1)(p-2)}p^{(p-2)(p+1)/2}\,.
\end{equation}
Calculation of $\dis(f_1)$:
\begin{align*}
 \dis(f_1) & = p^{(p-3)(p+1)/4}\prod_{0\leq i<j\leq (p-3)/2}\left(\omfroot{g^{2i}}-\omfroot{g^{2j}}\right)^2 \\
         & = p^{(p+1)(p-3)/4}\left(\frac{-1}{p}\right)^{(p-3)(p-1)/8}
\prod_{0\leq i<j\leq (p-3)/2}\left(\frac{2(\z^{g^{2i}}-\z^{g^{2j}})}{(1-\z^{g^{2i}})(1-\z^{g^{2j}})}\right)^2 \\
     & = -2^{(p-1)(p-3)/4}p^{(p-3)(p+3)/8}
  \prod_{0\leq i<j\leq (p-3)/2}\left(\frac{\z^{g^{2i}}-\z^{g^{2j}}}{(1-\z^{g^{2i}})(1-\z^{g^{2j}})}\right)^2\,.
\end{align*}
The right-most product in the last equality is a unit times $\la^{-(p-1)(p-3)/4}$, therefore, 
$\wal(\dis(f_1))=(p-1)(p-3)(p+3)/8 - (p-1)(p-3)/4= (p-1)(p-3)(p+1)/8$.
Since $\dis(f_1)$ is a rational integer, it follows that
\begin{equation} \label{discr f1}
 \dis(f_1)=\pm 2^{(p-1)(p-3)/4}p^{(p-3)(p+1)/8}\,.
\end{equation}
Now the relations (\ref{disc and res}), (\ref{discr f}) and (\ref{discr f1}) imply the validity
of the relation (\ref{resultant f1,f2}).
\begin{quote}
{\em Third Claim}: Among the integers $f_1(a),f_2(a)$ one is equal to $\pm 1$ and the other is equal 
to $\pm q^n$.
\end{quote}
Proof: By Bezout's identity, there exist polynomials $h_1(x),h_2(x)\in\Z[x]$ (both of degree $< (p-1)/2$) 
such that 
\[
 h_1(x)f_1(x)+h_2(x)f_2(x)=\dis(f_1,f_2)= \pm 2^{(p-1)^2/4} p^{(p^2-1)/8}\,.
\]
We make the substitution $x \leftarrow a$ in this equality.
By $f_1(a)f_2(a)=f(a)=\pm q^n$ and the fact that $(q,2p)=1$ (cf. beginning of the proof of Proposition 2.1), 
it follows that exactly one $f_i(a)$ is equal to $\pm 1$ and the other is equal to $\pm q^n$.
\begin{quote}
{\em Fourth Claim}: If $f_i(a)=\pm 1$ for some $i\in\{1,2\}$, then $a=0$.
\end{quote}
Proof: Let us put $f_1(x)=c_0+c_1x+\cdots +c_r x^r$, where $r=(p-1)/2$. We already know that $c_0=\pm 1$ and
$c_i\in\Z$ for all $i$. By the very definition of the polynomial $f_1$, its roots are 
\[
\xi_j =  \frac{1+\zeta^j}{\om (1-\zeta^j)}\,,\quad j\in S\,,
\]
where $S$ is a complete set of quadratic residues $\bmod\,p$.
\\
Since $\Phi(-1)=1$ the numerator $1+\z^j$ is a unit and it follows easily that
\[
\wal(\xi_j)= - (p+1)/2\,.
\]
Thus, for $k>0$,
\[
\wal(c_k) \ge -\frac{p+1}{2}\frac{p-1-2k}{2}+\frac{(p-1)(p+1)}{4} =\frac{k(p+1)}{2} \,.
\]
Then, for $k\geq 1$, we have $\val(c_k)=\wal(c_k)/(p-1)>0$ and therefore 
\begin{equation} \label{p divides all c}
 p \mid c_k\,,\quad k=1,\ldots,(p-1)/2\,.
\end{equation}
Moreover, for $k\geq 2$, we have $\val(c_k)\geq (p+1)/(p-1)> 1$, hence $\val(c_k)\geq 2$ and we can
write therefore
\begin{equation} \label{p-adic expansion}
f_1(x) \equiv c_0+c_1 x \pmod{p^2 x^2}\,.
\end{equation}
Next, we prove that
\begin{equation} \label{val c1}
 \val(c_1)=1\,.
\end{equation}
First, note that $f_2(x)=f_1(-x)$, which results from the fact that the polynomial $f_1$ is of odd degree,
and the polynomials $f_1$ and $f_2$ have opposite leading coefficients and roots 
(cf.~just before the relation (\ref{disc and res})). These observations imply that
$f_2(x)=c_0-c_1x+c_2x^2-c_3x^3+\cdots $. 
\\
Another observation is that $c_1\neq 0$. Indeed, since $f(x)=f_1(x)f_2(x)$, the coefficient of $x^2$ in
$f(x)$ is $2c_0c_2 +c_1^2$. On the other hand, by the initial definition of $f(x)$, the coefficient of $x^2$
is  $-p\bincf{p}{2}$, which is odd, because $p\equiv 3\pmod 4$. Therefore, $c_1$ is odd; in particular, it is 
non-zero.
\\
A third fact --which is a bit more than an observation-- is that 
\begin{equation} \label{ub for c1}
 |c_1| < p^2\,.
\end{equation}
If we prove this, then, in combination with the relation (\ref{p divides all c}) and the fact that $c_1\neq 0$, 
we will conclude that $\val(c_1)=1$.
\\ 
Proof of (\ref{ub for c1}): Let $g_1(x)$ be, as before, the reciprocal of the polynomial $f_1(x)$. 
Then $c_1$ is equal,
up to sign, with the sum of the roots of $g_1(x)$. But the roots of $g_1(x)$ are the reciprocals of the roots 
of $f_1(x)$, i.e. they are equal to $\xi_j^{-1}$, $i=1,\ldots, (p-1)/2$. 
Therefore (remember that $S$ is a complete set of residues $\bmod\,p$),
\[
 |c_1|\leq \sqrt{p}\sum_{j\in S}\left|\frac{1-\z^j}{1+\z^j}\right|
              =\sqrt{p}\sum_{j\in S}\left|\tan\frac{\pi j}{p}\right| =
          \sqrt{p}\sum_{j\in S}\frac{1}{\left|\tan\frac{\pi(p-2j)}{2p}\right|}\,.
\]
Since $\left|\frac{\pi(p-2j)}{2p}\right|<\frac{\pi}{2}$, it follows that 
$\left|\tan\frac{\pi(p-2j)}{2p}\right|>\frac{\pi}{2p}|p-2j|$, hence,
\[
 |c_1| <\frac{2p\sqrt{p}}{\pi}\sum_{j \in S} \frac{1}{|p-2j|}\,.
\]
Note that, as $j$ runs through the set $S$, the numbers $|p-2j|$ are distinct $\bmod\,p$, for, if
$|p-2j_1|\equiv |p-2j_2| \pmod p$ with $j_1,j_2\in S$ and $j_1\neq j_2$, then, necessarily, $j_2 = -j_1$,
which implies that $-1$ is a quadratic residue $\bmod\,p$, a contradiction. Therefore, the set 
$\{|p-2j|: j\in S\}$ is a subset of $\{1,\ldots,p-1\}$ with cardinality $(p-1)/2$. It is clear, therefore, that
\[
 \sum_{j \in S} \frac{1}{|p-2j|} \leq \sum_{k=1}^{(p-1)/2}\frac{1}{k} <\frac{3}{2}+\log\frac{p-1}{4}\,,
\]
from which we obtain
\[
 |c_1| <\frac{2p\sqrt{p}}{\pi}\left(\frac{3}{2}+\log\frac{p-1}{4}\right)\,.
\]
This upper bound for $|c_1|$ clearly implies $|c_1|<p^2$, as claimed.

{\em Final step of the proof of Theorem \ref{th p equiv 3 mod 8}}: By our third claim above, 
$f_1(a)$ or $f_2(a)$ must be $\pm 1$.
Since $f_2(a)=-f_1(-a)$, we may suppose that $f_1(a)=\pm 1$, i.e. $c_0+c_1a+c_2a^2+\cdots = \pm 1$.
Remember that $c_0=\pm 1$, therefore, $c_0+c_1a+c_2a^2+\cdots = \pm c_0$. The $-$ sign implies
$\pm 2+c_1a+c_2a^2+\cdots = 0$, clearly impossible, in view of (\ref{p divides all c}). The $+$ sign implies
$0=c_1a +c_2a^2 +\cdots $. If $a\neq 0$, then, taking also into account (\ref{p-adic expansion}), we
obtain $0=c_1\pmod{p^2}$ which contradicts (\ref{val c1}). 
This forces $a=0$ and then, by (\ref{q^n}), $q^n=\pm 2$, a contradiction.
\proofend

\section{The equation $5x^2-4=y^n$} \label{section general exp eq}
The third relation (\ref{shape of sols}), written as $q^n=5(20v^2-1)^2-4$, naturally leads to the study of
the more general equation 
\begin{equation} \label{general exp eq}
 5x^2 -4=y^n\,,\quad \textcolor{blue}{\hbox{$y>1,\,n>2$} } 
\end{equation}
in the integer unknowns $x,y,n$, where $x$ and $y$ have not, of course, the same meaning as the $x,y$ 
in equation (\ref{Abu eq}).

First, let $n$ be even. It is well-known that the positive integer solutions of $5X^2-4=Y^2$ are given 
by $X=F_{2k+1}$, $Y=L_{2k+1}$ for $k>0$, where $F$ denotes the Fibonacci and $L$ the Lucas sequence; notice that
$k=0$ gives $y=1$ which is excluded. Since it is known that the only Lucas number which is a pure power is 
$L_3=4$ (\cite{BMS}, Theorem 2), it follows that the only solution $(x,y,n)$ of the equation 
(\ref{general exp eq}) with even $n$ is $(2,2,4)$. 

From now on we suppose that $n$ is odd.
\\
If $n=3$, then $(25x)^2-500=(5y)^3$. It is well known that the only integral solutions $(X,Y)$ 
to $Y^2=X^3+500$ are $(X,Y)=(5,\pm 25)$, corresponding to $y=1$, which has been excluded by hypothesis. 
Hence we may assume that $\gcd(n,3)=1$, in particular $n\ge 5$.
\\
If $x$ is even then $\,5x^2-4 \equiv -4, \,16,\, 12 \pmod{32}$
implying $n\le 4$, which has already been excluded. Hence $x$ and $y$ are odd.

Now we work in the field $K=\Q(\th)$, where $\th=\sqrt{5}$.
From now on and until the end of the paper we view $K$ as embedded into the real numbers with 
$\th\mapsto \sqrt{5}=2.2360679\ldots$. 
The ring of integers in $K$ is $\bI=\{(x+y\th)/2\,;\, x, y \in \Z\ {\rm with}\ x\equiv y \pmod{2} \}$,
$\vep=(1+\sqrt 5)/2$ is the fundamental unit. In $K$ unique factorization holds. 
Throughout this section, for $\al\in K$, $\al'$ will always denote
the algebraic conjugate of $\al$. 
We factorize the equation (\ref{general exp eq}) over the field $K$
\begin{equation} \label{eq factorized over K}
5x^2-4=(x\th-2)(x\th+2)\,.
\end{equation}
If $p$ is a (rational) prime divisor of $y$, then $5x^2-4=y^n$ implies that $p$ is odd
and, clearly, $5$ is a quadratic residue $\bmod\,p$. It follows that $p$ splits in $K$ 
and $p\equiv \pm 1 \pmod 5$.
Therefore $y$ factorizes in $\bI$ as $y=\pi \pi'$, where we can choose $\pi>0$ (then $\pi'$ is also positive). 
Notice also that $y^n\equiv 1 \pmod 5$, hence $y\equiv 1 \pmod {10}$ (remember that $y$ and $n$ are odd) 
and $y\ge 11$.

Without loss of generality, we assume that $x$ is positive.
Since $x$ is odd, $x\th+2$ and $x\th-2$ are coprime with $x\th > y^{n/2}\geq 11^{5/2}$. Hence, there 
exists $k\in\Z$ such that $x\th+2=\vep^k\pi^n$. 
Writing $k=\ell n+k_1$ with $-(n-1)/2\leq k_1\leq (n-1)/2$, we have 
$x\th+2=\vep^{k_1}(\vep^{\ell}\pi)^n=\vep^{k_1}\pi_1^n$, where $\pi_1=\vep^{\ell}\pi$. The conjugate relation 
is $-x\th+2=\vep'^{k_1}\pi_1'^n$ and summing the two relations we get
\begin{equation} \label{eps and pi}
 \vep^{k_1}\pi_1^n + \vep'^{k_1}\pi_1'^n = 4\,.
\end{equation}
We have $\pi_1=u+v\th$ or $\pi_1=(u+v\th)/2$, where $u,v\in\Z$ are unknown and in the second case $uv$ is odd.
Then, for fixed $n$ and $k_1$ we obtain from (\ref{eps and pi})
\begin{equation} \label{Thue eq}
 T_{k_1}(u,v):= \vep^{k_1}(u+v\th)^n + \vep'^{k_1}(u-v\th)^n = \mbox{$4$ or $2^{n+2}$,}
\end{equation}
where, in the second case, $uv$ is odd.
Note that the left-hand side of (\ref{Thue eq}) is a homogeneous polynomial in $\Z[u,v]$ of degree $n$,
hence the relation (\ref{Thue eq}) implies a Thue equation. Since 
$T_{-k_1}(u,v)=(-1)^{k_1}T_{k_1}(u,-v)$, it suffices to consider the Thue equations
$T_{k_1}(u,v)=\pm 4$ and $T_{k_1}(u,v)=\pm 2^{n+2}$ with $k_1=0,1,\ldots,(n-1)/2$, where, in the second 
equation, $uv$ is odd. Moreover, since the degree of the form $T_{k_1}$ is odd, we can ignore the minus
sign in the right-hand sides. 
Using the above Thue equations we will prove that there are no solutions $(x,y,n)$ to (\ref{general exp eq})
with $n\in\{5,7,11,13\}$. Actually, we will show that for these values of $n$ the Thue equations 
$T_{k_1}(u,v)=2^{n+2}$ with $uv$ odd, and $T_{k_1}(u,v)=4$ are impossible for all $k_1=0,1,\ldots,(n-1)/2$. 
For every $n$ as above, the method is practically the same.
However, as one can guess, the case $n=13$ is somewhat more complicated; so we briefly expose this case
in order to illustrate how we work. Numerous Thue equations of degree $n$ arise.
A practical method for the solution of such equations has been developed since long by 
Tzanakis and de Weger \cite{TdW1} which later was improved by Bilu and Hanrot \cite{BH} and implemented in 
\textsc{Pari} ($\mathtt{http://pari.math.u-bordeaux.fr/}$ and 
\textsc{Magma} \cite{Bosma}, \cite{magma-handbook}. 
We use either of these packages to solve the Thue equations that arise.

We assume now that $n=13$ and we consider all $k_1$'s in $\{0,1,\ldots,6\}$.\\
$k_1=0$: Since $T_0(u,v)$ is reducible, our equations are treated by elementary means; no solutions arise.
\\
$k_1=1$: Both equations $T_1(u,v)=4, 2^{15}$ are easily solved.
\\
$k_1=2$: The congruences $T_2(u,v)\equiv 4, 2^{15}\pmod{13^2}$ are impossible.
\\
$k_1=3$: The equation $T_3(u,v)=2^{15}$ with $uv$ odd implies solvability of the congruence 
$T_3(x,1)\equiv 0\pmod{2^{14}}$. But, as it is easily checked, this congruence has no solutions. The equation
$T_3(u,v)=4$ remains. Since $T_3(u,v)=4u^{13}+\cdots$, we multiply by $2^{11}$ and we obtain a Thue equation
$u'^{13}+65u'^{12}v+\cdots + 320000000v^{13}=2^{13}$, whose only solution is $(u',v)=(2,0)$ which we
obviously reject.
\\
$k_1=4$: Now, $T_4(u,v)=7u^{13}+\cdots + 234375v^{13}$. On multiplying by $7^{12}$ we obtain monic Thue
equations with right-hand sides $4\cdot 7^{12}$ and $2^{15}7^{12}$. No solutions are returned.
\\
$k_1=5$: Similarly to the case $k_1=2$, both congruences $T_5(u,v)\equiv 4, 2^{15}\pmod{13^2}$ 
are impossible.
\\$k_1=6$: All coefficients of $T_6(u,v)$ are even and 
$\frac{1}{2}T_6(u,v)=9u^{13} +260u^{12}v+\cdots +312500v^{13}= T'_6(u,v)$, say. We thus have the Thue
equations $T'_6(u,v)=2^{14}$ with $uv$ odd, and $T'_6(u,v)=2^2$. The first equation implies solvability
of the congruence $T'_6(x,1)\equiv 0\pmod{2^{13}}$ with $x$ odd, which is impossible. For the second equation
we are obliged to multiply by $3^{24}$ in order to obtain a monic Thue equation, as required by both 
\textsc{Pari} and \textsc{Magma}. The resulting equation is treated with some ``effort'' by \textsc{Magma} and
no solutions are returned. On the other hand, \textsc{Pari} after several hours was still ``struggling'', so
we gave up.

The computational difficulties arising above, when $k_1=6$ show the limitation of the method and, indeed, for
$n=17$ the computational difficulties for the solution of the resulting Thue equations, at present, seem to be
insurmountable.

Summing up our results so far, we have the following theorem.
\begin{proposition} \label{prop small n}
There are no solutions $(x,y,n)$ to the equation (\ref{general exp eq}) with $n$ divisible by at least
one of the primes $3,5,7,11$ or $13$. 
The only solution $(x,y,n)$ to the equation (\ref{general exp eq}) with even $n$ is $(\pm 2,2,4)$.  
\end{proposition}

{\em Computing a first upper bound for $n$}. We now fix a solution $(x,y,n)$ of the equation (\ref{general exp eq}),
where, in view of Proposition \ref{prop small n}, we can assume that $n\geq 17$. Obviously, we can also assume that
$n$ is prime. Based on the few observations just after the equation (\ref{eq factorized over K}), but relaxing
the condition $x>0$, we see that there exists a set $\cP$ consisting of (unordered) sets $\{\pi,\pi'\}$ such 
that $\pi>0$, $\pi\pi'=y$ and, if $\{\pi_1,\pi_1'\}$ and $\{\pi_2,\pi_2'\}$ are distinct elements of $\cP$,
then $\pi_2,\pi_2'$ are non-associated to both $\pi_1,\pi_1'$. 
\\
We modify $\cP$ as follows: Let $\{\pi,\pi'\}\in\cP $. There exists precisely an $m\in\Z$ such that
$\vep^m\leq \sqrt{\vep y}/\pi < \vep^{m+1}$. The last relation is equivalent to 
$\vep^{2m-1}\pi^2\leq y<\vep^{2m+1}\pi^2$. On putting $\vep^m\pi=\pi_1$ we obtain 
\begin{equation} \label{ineq pi1}
 \frac{\pi_1}{\sqrt{\vep}} \leq \sqrt{y} <\pi_1\sqrt{\vep}\,,\quad
\mbox{or, equivalently, $\quad \displaystyle{\sqrt{\frac{y}{\vep}}<\pi_1\leq \sqrt{\vep y}}$.}
\end{equation}
Note that $\pi_1'=(-1)^m\vep^{-m}\pi'$, so that $\pi_1|\pi_1'|=y$. On multiplying the first 
relation (\ref{ineq pi1}) by $|\pi_1'|$ we get 
$\displaystyle{\frac{y}{\sqrt{\vep}}\leq |\pi_1'|\sqrt{y} < y\sqrt{\vep}}$, hence
$\displaystyle{\sqrt{\frac{y}{\vep}} \leq |\pi_1'| < \sqrt{y\vep}}$. The last relation combined with
the second relation (\ref{ineq pi1}) implies $\max\{\pi_1/|\pi_1'|\,,\,|\pi_1'|/|\pi_1|\} \leq\vep$ and,
certainly, the left-hand side of the last inequality is $>1$. We make the substitution
$\pi\leftarrow \pi_1$ or $\pi\leftarrow |\pi_1'|$ according as $\pi_1/|\pi_1'|$ is $>1$ or $<1$, respectively.
In this way, an ``adjusted'' set $\cP_1$ replaces the set $\cP$ containing elements $\{\pi,\pi'\}$ such that,
\begin{equation} \label{normalized pi}
 \pi>0\,, \quad \pi|\pi'|=y\,,\quad 1<\pi/|\pi'|\leq \vep\,.
\end{equation}
Now, in view of the relation (\ref{eq factorized over K}) and the fact that the two factors in the left-hand
side are relatively prime, we must have an ideal equation $\ideal{2+x\th}=\ideal{\pi}^n$ or
$\ideal{2+x\th}=\ideal{\pi'}^n$ for some $\{\pi,\pi'\}\in\cP_1$, and then $\ideal{-2+x\th}=\ideal{\pi'}^n$ or
$\ideal{-2+x\th}=\ideal{\pi}^n$, respectively. By choosing the appropriate sign for $x$ we may assume that
$\ideal{2+x\th}=\ideal{\pi}^n$ , from which it follows that
\begin{equation} \label{eq 2+xth}
 2+x\th=\si\vep^k\pi^n \quad\mbox{for some $k\in\Z$ and $\si\in\{-1,1\}$,}
\end{equation}
%
%
%
%
and
\begin{equation} \label{eq xth-2}
x\th - 2 =\frac{y^n}{x\th+2}=\frac{\pi^n|\pi'|^n}{\si\vep^k\pi^n}=\si\vep^{-k}|\pi'|^n \,.
\end{equation}
By (\ref{eq 2+xth}) and (\ref{eq xth-2}) we obtain
\begin{equation} \label{eq power-1}
 \vep^{2k}\left(\frac{\pi}{|\pi'|}\right)^n - 1 =\frac{\si\vep^k\pi^n}{\si\vep^{-k}|\pi'|^n} - 1 
  = \frac{x\th+2}{x\th-2}-1 =\frac{4}{x\th-2}\,.
\end{equation}
We have $5x^2=y^n+4$, from which $|x|\th >y^{n/2}$.

Now we put
\begin{equation} \label{def Lambda}
  \La = 2k\log\vep - n\log\frac{|\pi'|}{\pi} \,,
\end{equation}
so that $\displaystyle{\La=\log(\vep^{2k}\left(\frac{\pi}{|\pi'|}\right)^n)}$ and now, by (\ref{eq power-1}),
\[
 |e^{\La}-1| =\frac{4}{|x\th-2|}\leq \frac{4}{|x|\th-2}<\frac{4}{y^{n/2}-2}< \frac{4.0001}{y^{n/2}}\,.
\]
Notice that the right most side is less than $5.63\cdot 10^{-9}$ in view of the fact that $y\geq 11$ 
and $n\geq 17$. Therefore,
\[
 |\La| < 1.01 |e^{\La}-1| < \frac{4.0402}{y^{n/2}}\,.
\]
On the other hand, since the ideals $\ideal{\pi}$ and $\ideal{\pi'}$ are distinct, $\pi/|\pi'|$ is
not a unit and, consequently, $\La\neq 0$. Thus,
\begin{equation} \label{ubd Lambda}
 0< |\La| < \frac{4.0402}{y^{n/2}}\leq
 \textcolor{blue}{\frac{4.0402}{11^{7/2}}<5.683\cdot 10^{-9}}
\end{equation}
and 
\begin{equation} \label{ubd logLambda}
 \log|\La| < -\frac{n}{2}\log y + 1.3963\,.
\end{equation}
Now we compare $k$ and $n$ that appear in the linear form $\La$. 
\textcolor{blue}{
We already know that $n\geq 17$ and we show that $k<0$. Indeed, $k$ cannot be 
strictly positive for, otherwise, 
$|\Lambda|=2k\log\ep +n\log(\pi/|\pi'|) \geq 2\log\ep>0.9624$ which contradicts
\eqref{ubd Lambda}. Also, $k\neq 0$, because, if $k=0$, then, from  
\eqref{eq 2+xth} and \eqref{eq xth-2} we obtain
$\pi^n - \pi'^n=\pm 4$. This relation along with $\pi^n\pi'^n=y^n$ implies
that $\pi^n,\pi'^n$ are \emph{real} roots of $X^2\mp 4X+y^n$, therefore
$4-y^n\geq 0$ which contradicts the fact $y\geq 11$ and $n\geq 17$. In conclusion,
$k<0$ and  
\begin{equation} \label{Lambda for applying Laurent}
\Lambda = n\log(\pi/|\pi'|) -2|k|\log\ep .
\end{equation}
}
Further, by (\ref{ubd Lambda}), $|\La|<\textcolor{blue}{5.683\cdot 10^{-9}}$, therefore, in view also of (\ref{normalized pi}),
\[
 |k|=-k=-\frac{\La}{2\log\vep}+\frac{n}{2\log\vep}\log\frac{\pi}{|\pi'|}
<
\textcolor{blue}{
5.905\cdot 10^{-9}+\frac{n}{2\log\vep}\log\vep =5.905\cdot 10^{-9}+\frac{n}{2}\,,
}
\]
hence
\begin{equation} \label{compare k,n}
 |k| \leq \frac{n}{2} \,.
\end{equation}
Next, we consider the algebraic number $\eta:=\frac{\pi}{|\pi'|}$ appearing in $\La$. This number
is a root of the polynomial 
\[
 (\pi X-|\pi'|)(|\pi'|X-\pi)=yX^2-(\pi^2+\pi'^2)X+y = yX^2-(a^2\pm 2y)X+y \in
 \textcolor{blue}{\Z[X]}\,,
\]
where $a=\pi+\pi'\in\Z$. From this we easily see that 
\begin{equation} \label{height eta}
 \Ht(\eta) <\frac{1}{2}(\log y +\log\vep)\,.
\end{equation}
Finally, we are ready to calculate a first upper bound for $n$ using Corollary 2 of \cite{Laurent}.
In view of the relations (\ref{height eta}) and (\ref{normalized pi}) it is easy to estimate the quantities
that are involved in that corollary. Choosing the parameter $m$ that appears in the corollary equal to 20, 
and taking into account that $\max(2|k|,|n|\}=n$ (cf.\ref{compare k,n})), we easily find that, if $n\geq 15100$,
then
\[
 \log|\La|\geq \textcolor{blue}{-78.8}\left(\log n +\log\frac{\log y +\log\vep+1}{\log y+\log\vep}+0.38\right)^2(\log y+\log\vep)\,.
\]
This, combined with (\ref{ubd logLambda}), gives
\begin{equation} \label{test corollary 2}
78.8\left(\log n +\log\frac{\log y +\log\vep+1}{\log y+\log\vep}+0.38\right)^2(\log y+\log\vep)
-\frac{n}{2}\log y +1.3963 > 0\,.
\end{equation}
Since $y\geq 11$, we easily check that the inequality (\ref{test corollary 2}) can hold only if 
\begin{equation} \label{first ub for n}
 n < 2.2\times 10^4\,.
\end{equation}

{\em Proving that solutions with ``small'' $y$ cannot exist}. 
Now we go back to our equation (\ref{general exp eq}) and we assume that $(x,y)$ is a positive solution.
It is easily checked that this positiveness restriction does not prevent us from obtaining again the 
relations (\ref{eq power-1}) and (\ref{ubd Lambda}). We write the last inequality in the following shape:
\[
 \left|\frac{\log\eta}{\log\vep}-\frac{2k}{n}\right| < \frac{4.0402}{n\log\vep\cdot y^{n/2}}\,,
\quad \eta:=\frac{\pi}{\pi'}\,.
\]
The right-hand side is, obviously, less than $1/(2n^2)$, which shows that $2k/n$ is a convergent to the
continued fraction expansion of $\log\eta/\log\vep$ and, moreover, the denominator of this convergent is less
than $10^5$, in view of (\ref{first ub for n}). Let $a_0,a_1,a_2,\ldots$ be the partial quotients and
$p_0/q_0,p_1/q_1,p_2/q_2,\ldots$ the convergents to that expansion. Let $h$ be the first subscript such that 
$q_h\geq 10^5$. Then, $2k/n=p_m/q_m$ for some $m\in\{0,\ldots,h-1\}$. We have now
\[
 \frac{1}{(a_{i+1}+2)q_i^2} < \left|\frac{\log\eta}{\log\vep}-\frac{p_i}{q_i}\right|\,,
\]
hence,
\[
\frac{1}{(a_{i+1}+2)n^2} < \left|\frac{\log\eta}{\log\vep}-\frac{2k}{n}\right| 
                                                <\frac{4.0402}{n\log\vep\cdot y^{n/2}}\,,
\]
from which it follows that 
\begin{equation} \label{lb for A}
 4.0402 (A+2)n >\log\vep\cdot y^{n/2}\,,\quad A:=\max\{a_0,a_1,\ldots,a_h\}\,.
\end{equation}
For every $y\equiv 1\pmod{10}$ with $y<3\cdot 10^9$ and for every $\pi$ as above (there are $2^m$ such $\eta$' s, 
where $m$ is the number of rational prime divisors of $y$), we compute $\eta$ and the continued fraction 
expansion of the real number $\log\eta/\log\vep$, and we check the validity of the relation (\ref{lb for A}).
These computations can be performed with the routines of either \textsc{Pari} or  \textsc{Magma}. We stress the fact
that an ordinary precision is sufficient since the denominators of the checked convergents have at most
10 decimal digits. The whole task took around 30 hours of computations with \textsc{Pari} in a usual PC; with \textsc{Magma} it would take more time.
It turns out that, except possibly if $n\leq 11$, this relation is {\em not} satisfied. But we already know
that $n\geq 17$, hence we conclude:
\begin{quote}
 No solutions $(x,y,n)$ to (\ref{general exp eq}) exist with $n\geq 17$ and $y\leq 3\cdot 10^9$.
\end{quote}
\textcolor{blue}{
{\em Obtaining a smaller upper bound for $n$.\footnote{Updated on October 8, 
2020.\label{update foot}}}
In the published paper\footnote{J.~Number Theory {\bf 131} (2011) 1575--1596.}
we accomplish this by using Theorem 1 of \cite{Laurent}.
However, as A.~Koutsianas pointed out to us, our choice of the parameter $a_2$ in
that Theorem is incorrect; see Remark 1 below. 
As a consequence, here we revise and correct this paragraph of the paper, at the cost 
of obtaining a worse upper bound for $n$ ($n\leq 1153$ instead of $n\leq 811$). 
We thank A.~Koutsianas for his pointing out this mistake. 
} 

\textcolor{blue}{ 
We consider our linear form $\Lambda=n\log\eta-2|k|\log\ep$, where 
$\eta=\pi/|\pi'|$ (cf.~\eqref{Lambda for applying Laurent}). 
Now, we know that  $y> y_0:=3\cdot 10^9$ (this is very important!) and we apply 
M.~Laurent's Theorem 2 of \cite{Laurent} to $\Lambda$. 
In the notation of that theorem, 
$\alpha_1=\ep$, $\alpha_2=\eta$, $b_1=2k$, $b_2=n$.
We keep going with the notation of \cite[Theorem 2]{Laurent}:
We chose $\rho=1/2$ and $\mu=16$, 
$a_1=(\rho+1)\log\ep$, $a_2=(\rho+1)\log\ep+2\log (y_0)$ (for the choice of $a_2$ we
make use of \eqref{height eta}) and 
$h=\max\lbrace 2\left(\log\left(\left(\dfrac{1}{a_1}+\dfrac{1}{a_2}\right)n\right)
         +2\log\la+1.75\right)+0.06,\,\la\rbrace$,  
where (following the theorem) $\la=\sigma\log\rho$ and $\sigma=(1+2\mu-\mu^2)/2$.
Laurent's theorem implies a lower bound, say $-B(n)$, for $\log|\Lambda|$, where 
$B(n)$ is an \emph{explicit} positive function of $n$ with $B(n)=O(\log n)$.
This, combined with \eqref{ubd logLambda} gives $B(n)- n/2\log y + 1.3963>0$ which
is impossible if $n$ is ``sufficiently large''. Specifically, our computations showed
that the prime $n$ must not exceed $1153$. Thus we have proved the following:
\begin{theorem} \label{th gen exp eq}
Any integer solution of the equation $5x^2-4=y^n$ with $y>1$ and $n$ an odd prime, satisfies: 
\emph{(i)}  $17\leq n\leq 1153$  
and 
\emph{(ii)} $ y>3\cdot 10^9$. 
\end{theorem}
}

\textcolor{blue}{
{\bf Remarks} (October 8, 2020)
(1)
In the published version of the paper (J.~Number Theory {\bf 131} (2011), 1575--1596),
in order to reduce the initially obtained upper bound 
$n<2.2\cdot 10^4$ (see the paragraph ``Obtaining a smaller upper bound for $n$'',
p.p.~1587--1588), we first apply Laurent's Theorem 2, as above, without taking serious
care about the choice of the parameters $\rho$ and $\mu$, obtaining thus the reduced
bound $n\leq 6404$ (hence $n\leq 6397$ if we assume that $n$ is prime).
Then we turn to Laurent's Theorem 1 of the same paper \cite{Laurent} which we apply
repeatedly by choosing each time different values for the parameters 
$R_1,S_1,R_2,S_2$ of that theorem until we arrive to $n<820$ (hence to $n\leq 811$ 
if we assume that $n$ is prime). In our application of Laurent's Theorem 1, $a_1,a_2$
are same with those we chose above when we apply Laurent's Theorem 2, but now the
condition (2) of Theorem 2 implies that a positive function of $a_1,a_2$ 
(in the notation of that Theorem, this function is $gL(Ra_1+Sa_2)$ with $g,L,R,S>0$
positive constants) is bounded from above by a constant. This is absurd, because $a_2$ is a strictly increasing function of $y$ of which no upper bound is known.
Therefore, we cannot use Laurent's Theorem 1. We thank A.~Koutsianas who pointed
out to us this misuse.     
}

\textcolor{blue}{
(2)
The explicit function $n\mapsto B(n)- n/2\log y + 1.3963$ which 
we mentioned above was computed by a simple {\sc Maple} program and then we made experiments
with various values of $\rho$ and $\mu$ until we decide that $(\rho,\mu)=(1/2,16)$
implies the best upper bound for $n$. Independently, A.~Koutsianas wrote a 
\textsc{Magma} program --once again we thank him-- following a somewhat different 
strategy in order to compute an optimum upper bound for $n$, which ends-up with the 
same upper bound $n\leq 1153$.   
}

(3) In recent years, the so called ``modular approach'' to certain types of Diophantine equations 
--the Fermat equation being one of them-- turned out to be very succesful; see, for example, S.~Siksek's 
``The modular approach to Diophantine equations'', Chapter 15 in \cite{Cohen}. Our equation $5x^2-4=y^n$ 
resembles the Lebesgue-Nagell equation $x^2+D=y^n$, to which the modular method applies succesfully in
most cases; see \cite{BMSII}. However, as mentioned in \cite{BMSII}, the method is not succesful when 
$D=-a^2\pm 1$ because, in that case, there exists an obvious solution valid for every $n$. 
In the case of our equation we face a similar situation: the existence of the solution $(x,y)=(1,1)$ 
for every $n$ makes the application of the modular method ``hopeless'', according to S.~Siksek 
(private communication).
\section{The equation $5x^2-4=y^n$ when $y$ is prime} \label{section exp eq with y prime}
The main result of this section is the following
\begin{theorem}  \label{th with y=q}
Let $q$ be an odd prime. Then, for the solutions $(x,n)$ of the equation
\begin{equation}  \label{eq 5x^2-4=q^n}
 5x^2-4 = q^n\,,\quad x>0\,,\; n>0\,,n\neq 2
\end{equation}
the following are true: 

\emph{(i)}\; If $q\not\equiv 1\pmod{10}$, no solutions exist.

\emph{(ii)}\; If $q\leq 3\cdot 10^9$, no solutions with $n>2$ exist.

\emph{(iii)}\; If $(q+4)/5=\Box$, then $(x,n)=(\sqrt{\frac{q+4}{5}},1)$ is the only solution.

\emph{(iv)}\; If $(q+4)/5\neq\Box$, then at most one solution exists.

\emph{(v)}\; No solutions exist with $n>820$.

\emph{(vi)}\; No solutions exist with $n$ divisible by a prime from the set $\{2,3,5,7,11,13\}$.
\end{theorem}
The proof of this theorem follows from a straightforward combination of Theorem \ref{th gen exp eq},
already proved in Section \ref{section general exp eq} and Proposition \ref{at most one solution}, below. 
Therefore, the present section is essentially devoted to the proof of this proposition.

As noted in the beginning of Section \ref{section general exp eq}, the third relation (\ref{shape of sols}), 
written as $5(20v^2-1)^2-4=q^n$, led us to the more general equation (\ref{general exp eq}) for which  
Theorem \ref{th gen exp eq} holds. In this theorem, $y$ is general and not necessarily prime as the
equation $5(20v^2-1)^2-4=q^n$ would suggest. In this section, however, we will add the extra restriction
that the unknown $y$ in the equation (\ref{general exp eq}) be a prime, say $y=q$, and we will prove the following
theorem.
\begin{proposition} \label{at most one solution}
If $q$ is an odd prime, then the equation
\begin{equation}  \label{eq 5x^2=q^n+4}
 5x^2=q^n+4\,,\quad \mbox{$x,n$ positive integers, $n$ odd}
\end{equation}
has at most one solution if $(q+4)/5$ is not a perfect square and exactly one solution, namely,
$(x,n)=(\sqrt{(q+4)/5},1)$ if $(q+4)/5$ is a perfect square.
\end{proposition}
{\em Remark}: It is easy to see that the relation (\ref{eq 5x^2=q^n+4}) implies $q\equiv 1\pmod{10}$.

\proof
The proof of Proposition \ref{at most one solution} will be completed in three steps.
\\[1mm]
{\em  Step 1: The gap between two solutions of (\ref{eq 5x^2=q^n+4}).}
This step consists in proving that, if two solutions $(x,n)$, and $(x^{\prime},n^{\prime})$ exist, with
$n^{\prime}>n$, then $n^{\prime}$ must be ``very large'' compared to $n$; see (\ref{lower bound n'}).
We need first the following result.
\begin{lemma}
\label{2 plus x sqrt 5} Let $x$ be a positive integer and assume that
\begin{equation}
2+x\sqrt{5}=\xi^{a}\,, \label{hypothesis lemma 1}
\end{equation}
where $\xi$ is an algebraic integer in $\Q(\sqrt{5})$ and $a$ is an
integer $>1$. Then, $\xi=\frac{1+\sqrt{5}}{2}$, $a=3,x=1$.
\end{lemma}

\vspace{-5mm}
{\em Proof of the lemma}. There are two possibilities for $\xi$:
(I) Either $\xi=\frac{b+c\sqrt{5}}{2}$ with $b,c$ odd integers,
or (II) $\xi=b+c\sqrt{5}$ with $b,c$ arbitrary integers.
\\
After expansion of the right-hand side of (\ref{hypothesis lemma 1}) we obtain
\[
b^{a} + 5c^{2}\binom{a}{2}b^{a-2}+\cdots=%
\begin{cases}
2^{a+1} & \mbox{in case (I)}\\
2 & \mbox{in case (II)}
\end{cases}\,.
\]
It follows from this that, if $a$ is even, then an odd power of 2 is a square
$\bmod{\,5}$ which is impossible. Therefore, $a$ is odd and
\[
b^{a} + 5c^{2}\binom{a}{2}b^{a-2}+\cdots+5^{\frac{a-1}{2}}c^{a-1}\binom{a}{a-1}b =
\begin{cases}
2^{a+1} & \mbox{in case (I)}\\
2 & \mbox{in case (II)}
\end{cases}\,.
\]
Case (II) is impossible. Indeed, note that in the left-hand side all exponents
of $b$ are odd and all exponents of $c$ are even, hence $b>0$. Also, $b$
divides 2, hence $b=1$ or $2$. If $b=1$ then an obviously impossible
congruence $\bmod{\,5}$ results; and if $b=2$ then $2^{a}\leq2$ which implies
$a=1$, contrary to the hypothesis. \newline In case (I) we have, as before,
$b>0$, $b$ is odd and $b|2^{a+1}$. Hence, $b=1$ and we have
\[
2^{a+1}=1+5\binom{a}{2}c^{2}+\cdots+5^{\frac{a-1}{2}}\binom{a}{a-1}c^{a-1}
\geq2^{a+1}\,,
\]
where the last inequality is strict for every $c$, if $a\geq 5$ and for every
$c$ with $|c|>1$ when $a=3$. Thus, to avoid the contradiction we must conclude
that $a=3$ and $|c|=1$ from which it easily follows that $c=1$ and
$\xi=(1+\sqrt{5})/2$. This completes the proof of Lemma \ref{2 plus x sqrt 5}.

We put ${\theta} =(1+\sqrt{5})/2,{\theta} ^{\prime}=(1-\sqrt{5})/2$. These are
the roots of the polynomial $x^{2}-x-1$ and ${\theta} $ is the fundamental
unit of the ring of integers of $\Q(\theta)$. In general, for any
${\al}\in\Q(\theta)$ we denote by ${\al}^{\prime}$ the
conjugate of ${\al}$ under the isomorphism ${\theta} \mapsto{\theta}^{\prime}$.

Assume now that $(x,n)$ is a solution to equation (\ref{eq 5x^2=q^n+4}). Then
$(2+x\sqrt{5})(2-x\sqrt{5})=-q^{n}$ and it is clear that the factors in the
left-hand side are relatively prime as algebraic integers of $\Q(\theta)$.
Also, every (rational) prime dividing $q$ factors into two
distinct prime ideals. It follows then that there exists an algebraic integer
${\si}$ with norm $\pm q$ such that the following ideal relation is true:
$(2+x\sqrt{5})=({\si})^{n}$. Then, for some $r\in \Z$ we have the
element equation $2+x\sqrt{5}=\pm{\theta}^{r}{\si}^{n}$ and since we can
assume without loss of generality that ${\si}>0$, we finally get
\[
2+x\sqrt{5}={\theta}^{r}{\si}^{n}\quad
\mbox{along with the conjugate relation $2-x\sqrt{5}=\th'^r\si'^n$.}
\]
Combining the last two relations we obtain
\[
0<{\de}:=
\left(\frac{{\theta}^{\prime}}{{\theta}}\right)^{r}\left(\frac{{\si}^{\prime}}{{\si}}\right)^{n}+1
=\frac{4}{{\theta}^{r}{\si}^{n}}=\frac{4}{2+x\sqrt{5}}<\frac{1}{2}\,.
\]
Then,
$\frac{1}{2}<1-{\de}=-(\frac{{\theta}^{\prime}}{{\theta}})^{r}(\frac{{\si}^{\prime}}{{\si}})^{n}<1$
and in view of the inequality
$|\log(1-x)|<|x|(1+|x|)$ (valid for $|x|<1/2$) we obtain
\begin{equation} \label{linear form in logs 1}%
\left\vert-r\log\left\vert \frac{{\theta}}{{\theta}^{\prime}}\right\vert
+n\log\left\vert \frac{{\si}^{\prime}}{{\si}}\right\vert \right\vert
<{\de}(1+{\de})\,,\quad{\de}=\frac{4}{2+x\sqrt{5}}
=\frac{4}{2+\sqrt{q^{n}+4}}<\frac{4}{q^{n/2}}\,.
\end{equation}
Now, let $(x^{\prime},n^{\prime})$ another solution to (\ref{eq 5x^2=q^n+4})
with $n^{\prime}>n$, $n^{\prime}$ odd and $x^{\prime}>0$ (hence, $x^{\prime}>x$).
Exactly as before we have a relation $2+x^{\prime}\sqrt{5}={\theta}^{r^{\prime}}{\si}^{n^{\prime}}$
for a convenient $r^{\prime}\in\Z$ and
\begin{equation}
\left\vert -r^{\prime}\log\left\vert \frac{{\theta}}{{\theta}^{\prime}}
\right\vert +n^{\prime}\log\left\vert \frac{{\si}^{\prime}}{{\si}}
\right\vert \right\vert <{\de}^{\prime}(1+{\de}^{\prime})\,,\quad
{\de}^{\prime}=\frac{4}{2+\sqrt{q^{n^{\prime}}+4}}<\frac{4}{q^{n^{\prime}/2}}
<{\de}\,. \label{linear form in logs 2}
\end{equation}
Putting $u=\log|{\theta}/{\theta}^{\prime}|=\log((3+\sqrt{5})/2)$ and
eliminating the term $\log|{\si}^{\prime}/{\si}|$ from the inequalities
(\ref{linear form in logs 1}) and (\ref{linear form in logs 2}) we get
$|-rn^{\prime}+r^{\prime}n|u<n^{\prime}{\de}(1+{\de})+n{\de}^{\prime}(1+{\de}^{\prime})
<2n^{\prime}{\de}(1+{\de})$, i.e.
\begin{equation}
|-rn^{\prime}+r^{\prime}n|<\frac{2n^{\prime}}{u}{\de}(1+{\de})\,.
\label{linear form r and r'}%
\end{equation}
The left-hand side in (\ref{linear form r and r'}) is non-zero. Indeed, in the
opposite case we would have $\frac{r}{n}=\frac{r^{\prime}}{n^{\prime}}=\,\mbox{(say)\,}\frac{r_{1}}{n_{1}}$
with $(r_{1},n_{1})=1$. Then, $r=ar_{1},n=an_{1}$, $r^{\prime}=br_{1},n^{\prime}=bn_{1}$ for some positive
odd integers $a,b$ with $a<b$ and, moreover,
$2+x^{\prime}\sqrt{5}=({\theta}^{r_{1}}{\si}^{n_{1}})^{b}$.
By Lemma \ref{2 plus x sqrt 5} we conclude that $x^{\prime}=1$, contrary to the fact that
$x^{\prime}>x>1$.
\\
We conclude therefore that the left-hand side of (\ref{linear form r and r'}) is $\geq1$,
from which it follows that
\begin{equation}
n^{\prime}>\frac{u}{2}{\de}^{-1}(1+{\de})^{-1}, \label{lower bound n'}
\end{equation}
which shows that, the larger solution $n^{\prime}$ is \textquotedblleft far
away\textquotedblright\ from the smaller solution $n$; specifically, it is of
the size of $q^{n/2}$. This fact will play an important role below.
\\[1mm]
{\em  Step 2: Application of Hypergeometric Polynomials.}
At this second step we adapt to our equation the method of F.~Beukers in \cite{Beu1} and \cite{Beu2}. 
As a result we prove Lemma \ref{technical lemma} below, after which the final step for the proof of 
Proposition \ref{at most one solution} is not difficult. 
In that method one uses as a tool the hypergeometric polynomials, the properties of which we remind 
immediately below.
\\
Given the real numbers ${\al},{\be},{\gamma}$
where ${\gamma}$ is not zero or a negative integer, we define the
\emph{hypergeometric function} (with parameters ${\al},{\be},{\gamma}$)
\[
F({\al},{\be},{\gamma},z)=1+\frac{{\al}{\be}}{{\gamma}}z
+\sum_{k=2}^{\infty}\frac{1}{k!}\frac{{\al}({\al}+1)
\cdots({\al}+k-1){\be}({\be}+1)\cdots({\be}+k-1)}{{\gamma}({\gamma}+1)
\cdots({\gamma}+k-1)}z^{k}
\]
which converges for every complex number $z$ with $|z|<1$ and, in case that
${\gamma}>{\al}+{\be}$, it also converges for $z=1$. Let $n_{2}>n_{1}>0$
be integers. Put $n=n_{1}+n_{2}$ and define
\[
G(z)=F(-n_{2}-1/2,-n_{1},-n,z)\,,\quad H(z)=F(-n_{1}+1/2,-n_{2},-n,z)\,.
\]
By the definition of $G$ it is easy to see that
\[
G(z)=\sum_{k=0}^{n_{1}}\binom{n_{2}+1/2}{k}\binom{n_{1}}{k}\binom{n}{k}^{-1}(-z)^{k}\,,
\]
which, in particular, shows that, for any real number $z<0$, $G(z)$ is
positive. We will use the following properties:

\vspace{-5mm}
\begin{enumerate}
\item \label{prop 1} $G(z)$ and $H(z)$ are polynomials in $z$ of degrees
$n_{1} $ and $n_{2}$, respectively. Moreover, the polynomials $\binom{n}{n_{1}}G(4z) $ and
$\binom{n}{n_{1}}H(4z)$ have integer coefficients.
\item \label{prop 2} $|G(z)-(1-z)^{1/2}H(z)|<G(1)|z|^{n+1}$ for $|z|<1$.
\item \label{prop 3} $G(1)<G(z)<G(0)<1$ for $0<z<1$.
\item \label{prop 4} If $G^{*}(z)$ is the polynomial resulting from $G(z)$
when $n_{1},n_{2}$ are respectively replaced by $n_{1}+1,n_{2}+1$ and
$H^{*}(z)$ is defined analogously, then
\[
G^{*}(z)H(z)-G(z)H^{*}(z)=cz^{n+1}
\]
for some non-zero constant $c$.
\item \label{prop 5} $\displaystyle{|G(z)| < \left(  1+\frac{|z|}{2}\right)^{n_{2}+1}} $ for any $z$.
\end{enumerate}

\vspace{-4mm}
For the proof of the first four properties see Lemmas 1,2,3 and 4 in
\cite{Beu1}. For the proof of the fifth property see relation (1.10), page 226
of \cite{TzW}.
Now we are in a position to prove the main result of this step.

\begin{lemma} \label{technical lemma}
Let $(x,n)$ be a solution to (\ref{eq 5x^2=q^n+4}),
where, as always, $x>0$ and $n\geq1$ is odd; we assume, moreover, that
$q^{n}>600$. Let $r,s$ be positive integers such that $q^{n}\geq2^{6+4s/r}$
and define the positive real number $\nu$ by means of the relation
\[
q^{n\nu}=2.007\times(4.03q^{n})^{r/s}\,.
\]
Finally, let $N=q^{n^{\prime}}$ where $n^{\prime}>n$ and let $y$ be any
integer. Then,
\[
\left\vert \frac{y\sqrt{5}}{N^{1/2}}-1\right\vert >\frac{0.27}{q^{n(3+\nu/2)}}q^{n/s}N^{-(1+\nu)/2}\,.
\]
\end{lemma}
{\em Proof of the lemma}. Let $n_{2}>n_{1}$ be positive integers which will be specified later
and $m=n_{1}+n_{2}$. Put $z=-q^{-n}$. Then, $|4z|<1$ so that $G(4z)$ and
$H(4z)$ are meaningful. By properties \ref{prop 2} and \ref{prop 3} of the
polynomials $G$ and $H$ we have
\[
\left\vert G(4z)-H(4z)(1+\frac{4}{q^{n}})^{1/2}\right\vert <G(1)\left(\frac{4}{q^{n}}\right)^{m+1}
<\left(\frac{4}{q^{n}}\right)^{m+1}\,,
\]
hence
\begin{equation} \label{inequality with G,H}%
\left\vert \binom{m}{n_{1}}G(4z)-\binom{m}{n_{1}}H(4z)\frac{x\sqrt{5}}{q^{n/2}}\right\vert
<\binom{m}{n_{1}}\left(\frac{4}{q^{n}}\right)^{m+1}\,.
\end{equation}
By property \ref{prop 1} and the fact that $G(x)>0$ for any negative real number $x$,
\[
\binom{m}{n_{1}}G(4z)=\frac{A}{q^{nn_{1}}}\quad\mbox{for some positive $A\in\Z$}
\]
and similarly,
\[
\binom{m}{n_{1}}H(4z)=\frac{B}{q^{nn_{2}}}\quad\mbox{for some $B\in\Z$.}
\]
Then, (\ref{inequality with G,H}) implies
$\left\vert \dfrac{A}{q^{nn_{1}}}-\dfrac{Bx}{q^{nn_{2}}}\dfrac{\sqrt{5}}{q^{n/2}}\right\vert
< \dbinom{m}{n_{1}}\left(\dfrac{4}{q^{n}}\right)^{m+1}$, from which
\[
\left\vert 1-\frac{Bx\sqrt{5}}{Aq^{n(n_{2}-n_{1}+1/2)}}\right\vert
<2^{m-1}\frac{q^{nn_{1}}}{A}\frac{2^{2(m+1)}}{q^{n(m+1)}}
=\frac{2^{3m+1}}{Aq^{n(n_{2}+1)}}\,.
\]
Now, let us put ${\ep}=\left\vert \frac{y\sqrt{5}}{N^{1/2}}-1\right\vert$, so that,
from the above inequality we have
\begin{equation}
\left\vert \frac{y}{N^{1/2}}-\frac{Bx}{Aq^{n(n_{2}-n_{1}+1/2)}}\right\vert
<\frac{1}{\sqrt{5}}\left( {\ep}+\frac{2^{3m+1}}{Aq^{n(n_{2}+1)}}\right)\,. \label{estimate difference}%
\end{equation}
Let ${\lambda}=\lceil\frac{n^{\prime}-n}{2n}\rceil$. Then,
\begin{equation}
q^{n({\lambda}-1)}<\left(  \frac{N}{q^{n}}\right)^{1/2}\leq q^{n{\lambda}}\,.
\label{define la}%
\end{equation}
Now comes the moment to choose $n_{1},n_{2}$. First we choose $n_{1}$ to satisfy
\begin{equation}
\frac{r}{s}{\lambda}\leq n_{1}\leq\frac{r}{s}{\lambda}+\frac{2s-1}{s}\,.
\label{choose n1}%
\end{equation}
We must keep in mind that there are exactly two consecutive positive integers
in the interval $[r{\lambda}/s\,,\,(r{\lambda}+2s-1)/s]$; this is a simple
exercise. Choose now $n_{2}$ by setting $n_{2}=n_{1}+{\lambda}>n_{1}$ and
remember that $m=n_{1}+n_{2}=2n_{1}+{\lambda}$. Moreover, we will need below
that the left-hand side of (\ref{estimate difference}) be non-zero. In the
next lines we show that we can choose $n_{1}$ in such a way that this
requirement be satisfied. \newline Suppose that for the smaller integer
$n_{1}$ in the interval $[r{\lambda}/s\,,\,r{\lambda}+2s-1)/s]$ the left-hand
side of (\ref{estimate difference}) is zero. Then, we can repeat the above
process with $n_{1}^{\prime}:=n_{1}+1$ in place of $n_{1}$ ($n_{1}^{\prime}$
still belongs to this interval), $n_{2}^{\prime}:=n_{2}+1$ in place of $n_{2}$
and $m^{\prime}:=n_{1}^{\prime}+n_{2}^{\prime}$ in place of $m$, so that the
polynomials $G$ and $H$ will be replaced by $G^{\ast}$ and $H^{\ast}$
respectively, and the integers $A,B$ by some other integers, say, $A^{\ast},B^{\ast}$.
Then, we will obtain an inequality analogous to (\ref{estimate difference}), namely,
\[
\left\vert \frac{y}{N^{1/2}}-\frac{B^{\ast}x}{A^{\ast}q^{n(n_{2}^{\prime}-n_{1}^{\prime}+1/2)}}\right\vert
<\frac{1}{\sqrt{5}}\left({\ep}+\frac{2^{3m^{\prime}+1}}{A^{\ast}q^{n(n_{2}^{\prime}+1)}}\right)\,.
\]
If the left-hand side were again zero, then we would have $B/A=B^{\ast}/A^{\ast}$
(note that $n_{2}^{\prime}-n_{1}^{\prime}=n_{2}-n_{1}$), which
would easily imply that $z=-4/q^{n}$ is a zero of the function $G^{\ast}\cdot H-G\cdot H^{\ast}$
and this contradicts property \ref{prop 4} of the
polynomials $G,H$. We conclude therefore that for at least one integer $n_{1}$
satisfying (\ref{choose n1}), the left-hand side of (\ref{estimate difference})
is non-zero and from now on we assume that we have selected such an $n_{1}$.

We now rewrite the term $Aq^{n(n_{2}-n_{1}+1/2)}$ appearing in the left-hand
side of (\ref{estimate difference}). We first observe that (\ref{define la})
implies $q^{n^{\prime}/2}\leq q^{n({\lambda}+1/2)}$ which shows that
$q^{n(2{\lambda}+1)}=q^{n^{\prime}}q^{2\mu}$ for some non-negative integer
$\mu$. Consequently, on putting $q^{\mu}=A_{0}$ (a positive integer), we have
\[
Aq^{n(n_{2}-n_{1}+1/2)}=Aq^{n({\lambda}+1/2)}=Aq^{n^{\prime}/2}A_{0}
=A_{0}AN^{1/2}\,.
\]
Going back to (\ref{estimate difference}), we get
\begin{align*}
\frac{1}{\sqrt{5}}\left(  {\ep}+\frac{2^{3m+1}}{Aq^{n(n_{2}+1)}}\right)
&  >\left\vert \frac{y}{N^{1/2}}-\frac{Bx}{Aq^{n(n_{2}-n_{1}+1/2)}}\right\vert
=\left\vert \frac{y}{N^{1/2}}-\frac{Bx}{A_{0}AN^{1/2}}\right\vert
=\frac{|A_{0}Ay-Bx|}{A_{0}|A|N^{1/2}}\\
&  \geq\frac{1}{A_{0}|A|N^{1/2}}=\frac{1}{|A|q^{n({\lambda}+1/2)}}\,,
\end{align*}
from which
\begin{equation}
{\ep}|A|q^{n({\lambda}+1/2)}+2^{3m+1}q^{-n(n_{1}+1/2)}>\sqrt{5}\,.
\label{sum larger than 1}%
\end{equation}
We estimate separately the second summand in the left-hand side of
(\ref{sum larger than 1}). By the hypothesis on the lower bound of $q^{n}$ and
(\ref{define la}) we have
\[
\frac{2^{3m+1}}{q^{nn_{1}}}<\frac{2^{6n_{1}+3{\lambda}+1}}{q^{nn_{1}}}
\leq\frac{2^{6n_{1}+3{\lambda}+1}}{2^{(6+4s/r)n_{1}}}=2^{3{\lambda}+1-4sn_{1}/r}
\leq2^{3{\lambda}+1-4{\lambda}}=2^{1-{\lambda}}\leq 1\,,
\]
which shows that the second summand in the left-hand side of
(\ref{sum larger than 1}) is $\leq q^{-n/2}<600^{-1/2}$. This shows that the
first summand in the left-hand side of (\ref{sum larger than 1}) is larger
than $5^{1/2}-600^{-1/2}>2.195$. Then, remembering also how $A$ has been
defined and using property \ref{prop 5} of the polynomial $G$, we get
\begin{align*}
2.195  &  \leq{\ep}|A|q^{n({\lambda}+1/2)}={\ep}q^{n(n_{1}+{\lambda}+1/2)}\binom{m}{n_{1}}|G(-4/q^{n})| \\
&  <{\ep}q^{n(n_{1}+{\lambda}+1/2)}\cdot 2^{m-1}\left(1+\frac{2}{q^{n}}\right)^{n_{2}+1} \\
&  \leq{\ep}q^{n(n_{1}+{\lambda}+1/2)}\cdot 2^{m-1}\left(1+\frac{2}{q^{n}}\right)^{m}
={\ep}q^{n/2}q^{n(n_{1}+{\lambda})}\cdot\frac{1}{2}\left(2+\frac{4}{q^{n}}\right)^{m} \\
&  <\frac{{\ep}}{2}q^{n/2}q^{n(n_{1}+{\lambda})}\times2.007^{m}
\quad\mbox{(since $q^n>600$)} \\
&  =\frac{{\ep}}{2}q^{n/2}q^{n{\lambda}(1+n_{1}/{\lambda})}\times2.007^{{\lambda}(1+2n_{1}/{\lambda})} \\
&  \leq\frac{{\ep}}{2}q^{n/2}q^{n{\lambda}(\frac{r}{s}
+\frac{2s-1}{{\lambda}s}+1)}\times2.007^{{\lambda}(\frac{2r}{s}+\frac{4s-2}{{\lambda}s}+1)} \\
&  =\frac{{\ep}}{2}q^{n/2}q^{n(2s-1)/s}\times2.007^{(4s-2)/s}\cdot(q^{n(1+r/s)}\times2.007^{1+2r/s})^{{\lambda}}\\
&  <\frac{{\ep}}{2}q^{n/2}q^{n(2s-1)/s}\times2.007^{(4s-2)/s}\cdot(2.007\cdot(4.03q^{n})^{r/s}q^{n})^{{\lambda}}\\
&  =\frac{{\ep}}{2}q^{n/2}q^{n(2s-1)/s}\times2.007^{(4s-2)/s}\cdot q^{n(1+\nu){\lambda}}
\quad\mbox{(by the definition of $\nu$)}\\
&  <\frac{{\ep}}{2}q^{n/2}q^{n(2s-1)/s}\times2.007^{4}\cdot q^{n(1+\nu){\lambda}}\,,
\end{align*}
from which we immediately get
\[
0.27q^{n(-\frac{5}{2}+\frac{1}{s})}<{\ep}q^{n{\lambda}(1+\nu)}<{\ep}(Nq^{n})^{(\nu+1)/2}
\]
(the right-most inequality being true because $q^{n{\lambda}}<(Nq^{n})^{1/2}$
in view of (\ref{define la})), and hence the claimed lower bound for
${\ep}=|\frac{y\sqrt{5}}{N^{1/2}}-1|$. This completes the proof of Lemma \ref{technical lemma}.
\\[1mm]
{\em  Step 3: Completion of the proof of Proposition \ref{at most one solution}}.
Assume that, if $(q+4)/5$ is not a perfect square there exists a solution to equation (\ref{eq 5x^2=q^n+4})
and if $(q+4)/5$ is a perfect square there exists a solution to this equation 
besides the obvious one which results from the relation $5(\sqrt{(q+4)/5})^2=q+4$. 
Thus, in both cases, our assumptions in particular imply that there exists a solution $(x_0,n_0)$ with $n_0>1$,
hence, by Theorem \ref{th gen exp eq}, we must have $q>3\cdot 10^9$.
Let $(x,n)$ be the least solution to equation (\ref{eq 5x^2=q^n+4}).
In order to prove the theorem, we will assume that a larger solution $(x',n')$ to (\ref{eq 5x^2=q^n+4})
exists and we will arrive at a contradiction.
\\
We put $N=q^{n^{\prime}}$, so that $N+4=5x^{\prime2}$, from which we get
\[
\frac{4}{N}=\left(  \frac{\sqrt{5}x^{\prime}}{N^{1/2}}-1\right)
\left(\frac{\sqrt{5}x^{\prime}}{N^{1/2}}+1\right)
=\left(  \frac{\sqrt{5}x^{\prime}}{N^{1/2}}-1\right)
\left(  \frac{\sqrt{N+4}}{\sqrt{N}}+1\right)
> 2\left(\frac{\sqrt{5}x^{\prime}}{N^{1/2}}-1\right)\,,
\]
therefore
\[
0<\frac{\sqrt{5}x^{\prime}}{N^{1/2}}-1<\frac{2}{N}\,.
\]
We apply Lemma \ref{technical lemma} with $y=x^{\prime}$, $r=1,s=2$; then it
is easy to check that $\nu<0.7178$ and by the conclusion of the lemma and the
last displayed inequality we get
\[
2N^{-1}>\frac{\sqrt{5}x^{\prime}}{N^{1/2}}-1>0.27\times q^{-n(5+\nu)/2}N^{-(1+\nu)/2}\,,
\]
hence
\[
\frac{(1-\nu)}{2}n^{\prime}\log q<\log(7.408)+\frac{5+\nu}{2}n\log q<\frac{5.627+\nu}{2}n\log q\,,
\]
from which
\textcolor{blue}{
\begin{equation}  \label{bound n' by n}
n^{\prime}<\frac{5.627+\nu}{1-\nu}n\,. 
\end{equation}
}
On the other hand, recalling that $u=\log((3+\sqrt{5})/2)$ and
${\de=}\frac{4}{2+\sqrt{q^{n}+4}}$,  we have in view of (\ref{lower bound n'}),
\begin{align*}
n^{\prime}  &  >\frac{u}{2}{\de}^{-1}(1+{\de})^{-1}
=\frac{u}{2}\cdot\frac{2+\sqrt{q^{n}+4}}{4}\left(1+\frac{4}{2+\sqrt{q^{n}+4}}\right)^{-1}  \\
&  =\frac{u}{8}\cdot\frac{(2+\sqrt{q^{n}+4})^{2}}{6+\sqrt{q^{n}+4}}
> \frac{u}{8}\cdot\frac{(2+q^{n/2})^{2}}{6+q^{n/2}}
>0.12\times\frac{(2+q^{n/2})^{2}}{6+q^{n/2}}\,.
\end{align*}
Combining this lower bound for $n^{\prime}$ with 
\textcolor{blue}{(\ref{bound n' by n})}, we
get the following relation:
\[
0.12\times\frac{(2+q^{n/2})^{2}}{6+q^{n/2}}<\frac{5.627+\nu}{1-\nu}n\,.
\]
By the definition of $\nu$, $(q^{n})^{\nu}=2.007\times(4.03q^{n})^{1/2}$. Solving for $\nu$ and
substituting into the above inequality we obtain
\begin{equation} \label{impossible ineq}%
 n\frac{6.127n\log q+\ga}{0.5n\log q-\ga}>0.12\frac{(2+q^{n/2})^2}{6+q^{n/2}}\,,\quad
                        \ga=\log (2.007\cdot 4.03^{0.5})\,.
\end{equation}
However, in view of the large size of $q$ we easily check that (\ref{impossible ineq}) is {\em not}
satisfied and this contradiction proves that the solution $(x^{\prime},n^{\prime})$ cannot exist.
\proofend

\end{document}